\documentclass[12pt]{article}
\usepackage{amsmath}
\usepackage{amssymb}
\usepackage{amsthm}
\usepackage{graphicx}
\usepackage{verbatim}
\usepackage{subfig,url}

\theoremstyle{plain}
\newtheorem {thm}{Theorem}[section]

\newtheorem {lem}[thm]{Lemma}
\newtheorem {pr}[thm]{Proposition}
\newtheorem {cor}[thm]{Corollary}
\newtheorem {obs}[thm]{Observation}
\newtheorem{defn}[thm]{Definition}

\newtheorem{conj}{Conjecture}[section]
\newtheorem{question}{Question}[section]

\theoremstyle{definition}
\newtheorem{example}{Example}[section]

\theoremstyle{remark}
\newtheorem{remark}[thm]{Remark}

\newcommand{\Em}[1]{\textbf{#1}}

\DeclareMathOperator{\Arg}{Arg}

\font\elevenss=cmss11

\font\eightss=cmss8

\font\sixss=cmss8 at 6pt

\newfam\ssfam
\textfont\ssfam=\elevenss \scriptfont\ssfam=\eightss
 \scriptscriptfont\ssfam=\sixss

\def\ee{\epsilon}

\def\vv{{\bf v}}
\def\uu{{\bf u}}

\def\P{{\mathbb P}}
\def\E{{\mathbb E}}
\def\F{{\cal F}}

\def\ee{\epsilon}

\def\one{{\bf 1}}

\def\C{{\mathbb C}}
\def\Z{{\mathbb Z}}

\def\half{\frac{1}{2}}
\def\cE{{\cal E}}
\def\zeros{{\cal Z}}
\def\vbar{{\underline{V}}}
\def\Cox{\hfill \Box}
\def\G{{\cal G}}
\def\gt{\tilde{g}}
\def\htt{\tilde{h}}

\makeatletter \def\romenumi{ \def\theenumi{\roman{enumi}}
\def\p@enumi{\theenumi} \def\labelenumi{(\@roman\c@enumi)}}
\makeatother

\begin{document}

\begin{titlepage}
\begin{center}
{\large \bf The distribution of zeros of the derivative of a
random polynomial}
\end{center}
\vspace{5ex}
\begin{flushright}
Robin Pemantle 
\footnote{Research supported in part by National Science Foundation 
grant \# DMS 0905937.}$^,$\footnote{University of Pennsylvania, Department of
Mathematics, 209 S. 33rd Street, Philadelphia, PA 19104, USA} \\
Igor Rivin
\footnote{Research supported in part by National Science Foundation Grant 
\# DMR 835586t, the Institute for Advanced Study, the Berlin Mathematics 
Society, and Technische Universit\"at Berlin.}$^,$\footnote{Temple University, 
Department of Mathematics, 1805 North Broad St Philadelphia, PA 19122} \\
\end{flushright}
\begin{flushright}
{\tt pemantle@math.upenn.edu, rivin@math.temple,edu}
\end{flushright}

\vfill

	{\bf ABSTRACT:} \hfill \\[2ex]
	In this note we initiate the probabilistic study of the critical points 
	of polynomials of large degree with a given distribution of roots. Namely,
	let $f$ be a polynomial of degree $n$ whose zeros are chosen
	IID from a probability measure $\mu$ on $\C$.  We conjecture 
	that the zero set of $f'$ always converges in distribution to 
	$\mu$ as $n \to \infty$.  We prove this for measures with finite 
	one-dimensional energy.  When $\mu$ is uniform on the unit 
	circle this condition fails.  In this special case the zero
	set of $f'$ converges in distribution to that the IID Gaussian 
	random power series, a well known determinantal point process.
	\vfill



	\end{titlepage}

	\setcounter{equation}{0}
	\section{Introduction} \label{sec:intro}

	Since Gauss, there has been considerable interest in the location 
	of the \emph{critical points} (zeros of the derivative) of polynomials 
	whose zeros were known -- Gauss noted that these critical points were 
	points of equilibrium of the electrical field whose charges were placed 
	at the zeros of the polynomial, and this immediately leads to the proof 
	of the well-known Gauss-Lucas Theorem, which states that the critical 
	points of a polynomial $f$ lie in the convex hull of the zeros of $f$  
	(see, e.g.~\cite[Theorem~6,1]{marden1949}).  There are too many 
	refinements of this result to state.  A partial list (of which
	several have precisely the same title!) is as 
	follows:~\cite{MR700266,MR1063061,MR1400355,MR1473453,MR874708,MR788579,MR0308374,MR0239062,MR0239061,MR0021148,MR0019157,MR1506859,MR0133437,MR1452801,MR2159699,MR2063118,MR2722582,MR1856819}).  
	Among these, we mention two extensions that are easy to state.
	\begin{itemize}
	\item Jensen's theorem: if $p(z)$ has real coefficients, then the non-real critical points of $p$ lie in the union of the ``Jensen Disks'', which are disks one of whose diameters is the segment joining a pair of conjugate (non-real) roots of $p.$
	\item Marden's theorem: 
	Suppose the zeroes $z_1, z_2,$ and $z_3$ of a third-degree polynomial $p(z)$ are non-collinear. There is a unique ellipse inscribed in the triangle with vertices $z_1, z_2, z_3$ and tangent to the sides at their midpoints: the Steiner inellipse. The foci of that ellipse are the zeroes of the derivative $p^\prime(z).$
	\end{itemize}
	There has not been any \emph{probabilistic} study of critical 
points (despite the obvious statistical physics connection) from this
viewpoint. There has been a very extensive study of random 
polynomials (some of it quoted further down in this paper), but 
generally this has meant some distribution on the coefficients of 
the polynomial, and not its roots. Let us now define our problem:

	Let $\mu$ be a probability measure on the complex numbers.
	Let $\{ X_n : n \geq 0 \}$ be random variables on a probability
	space $(\Omega , \F , \P)$ that are IID with common distribution 
	$\mu$.  Let 
	$$f_n (z) := \prod_{j=1}^n (z - X_j)$$
	be the random polynomial whose roots are $X_1 , \ldots , X_n$.
	For any polynomial $f$ we let $\zeros (f)$ denote the empirical
	distribution of the roots of $f$, for example,
	$\zeros (f_n) = \frac{1}{n} \sum_{j=1}^n \delta_{X_j}$.

	The question we address in this paper is: 
	\begin{question} 
	\label{mainq}
	When are the zeros
	of $f_n'$ stochastically similar to the zeros of $f_n$?
	\end{question}
	 Some
	examples show why we expect this.  
	\begin{example}
	\label{eg1}
	Suppose $\mu$ concentrates
	on real numbers.  Then $f_n$ has all real zeros and the zeros of
	$f_n'$ interlace the zeros of $f_n$.  It is immediate from this
	that the empirical distribution of the zeros of $f_n'$ converges
	to $\mu$ as $n \to \infty$.  The same is true when $\mu$ is
concentrated on any affine line in the complex plane: interlacing 
holds and implies convergence of the zeros of $f_n'$ to $\mu$.  
Once the support of $\mu$ is not 
contained in an affine subspace, however, the best we can say 
geometrically about the roots of $f_n'$ is that they are contained 
in the convex hull of the roots of $f_n$; this is the 
Gauss-Lucas Theorem.
\end{example}

\begin{example}
\label{atomiceg}
Suppose the measure $\mu$ is
atomic.  If $\mu (a) = p > 0$ then the multiplicity of $a$ as
a zero of $f_n$ is $n (p + o(1))$.  The mulitplicity of $a$
as a zero of $f_n'$ is one less than the multplicity as a 
zero of $f_n$, hence also $n (p + o(1))$.  This is true for
each of the countably many atoms, whence it follows again
that the empirical distribution of the zeros of $f_n'$ converges
to $\mu$.  
\end{example}

Atomic measures are weakly dense in the space of
all measures.  Sufficient continuity of the roots of $f'$ 
with respect to the roots of $f$ would therefore imply that
the zeros of $f_n'$ always converge in distribution to $\mu$ 
as $n \to \infty$.  In fact we conjecture this to be true.

\begin{example}
Our first experimental example has the roots of $f$ uniformly distributed in the unit disk. In the figure, we sample $300$ points from the uniform distribution in the disk, and plot the critical points (see Figure~\ref{disksamp}).
\begin{figure}
\centering
\includegraphics[width=3in]{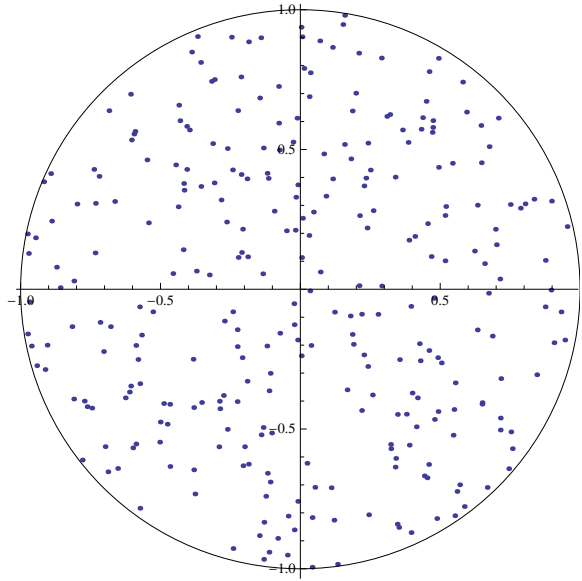}
\caption{Critical points of a polynomial whose roots are uniformly sampled inside the unit disk.}
\label{disksamp}
\end{figure}
The reader may or may not be convinced that the critical points are uniformly distributed.
\end{example}

\begin{example}
Our second example takes polynomials with roots uniformly distributed 
on the unit \emph{circle}, and computes the critical points. 
In Figure~\ref{circlesamp} we do this with a sample of size 300.
\begin{figure}
\centering
\includegraphics[width=3in]{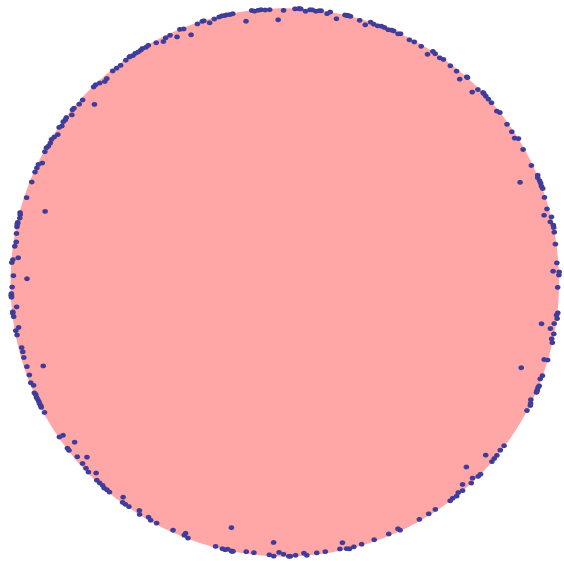}
\caption{Critical points of polynomial whose roots are uniformly sampled on the unit circle.}
\label{circlesamp}
\end{figure}
One sees that the convergence is rather quick.
\end{example}
\begin{remark}
The figures were produced with Mathematica. However, the reader wishing to try this at home should remember to do all the computations to very high precision, otherwise Mathematica (and Matlab and R) return garbage answers, since they all use a rather primitive method of computing zeros of polynomials.
\end{remark}
\begin{conj} \label{conj:1}
For any $\mu$, as $n \to \infty$, $\zeros (f')$ converges 
weakly to $\mu$.
\end{conj}

There may indeed be such a continuity argument, though the following
counterexample shows that one would at least need to rule out some
exceptional sets of low probability.  Suppose that $f(z) = z^n -1$.
As $n \to \infty$, the distribution of the roots of $f$ converge
weakly to the uniform distribution on the unit circle.  The roots
of $f_n'$ however are all concentrated at the origin.  If one moves
one of the $n$ roots of $f_n$ along the unit circle, until it meets
the next root, a distance of order $1/n$, then one root of $f_n'$
zooms from the origin out to the unit circle.  This shows that 
small perturbations in the roots of $f$ can lead to large perturbations
in the roots of $f'$.  It seems possible, though, that this is only
true for a ``small'' set of ``bad'' functions $f$.  

\subsection{A little history}
This circle of questions was first raised in discussions between 
one of us (IR) and the late Oded Schramm, when IR was visiting 
at Microsoft Research  for the auspicious week of 9/11/2001. 
Schramm and IR had some ideas on how to approach the questions, 
but were somewhat stuck.  There was always an intent to return 
to these questions, but Schramm's passing in September 2008 threw 
the plans into chaos. We (RP and IR) hope we can do justice to 
Oded's memory.  

These questions are reminscent of questions of the 
kind often raised by Herb Wilf, that sound simple but are not.  
This work was first presented at a conference in Herb's honor 
and we hope it serves as a fitting tribute to Herb as well.

\setcounter{equation}{0}
\section{Results and Notations}

Our goal in this paper is to prove cases of Conjecture~\ref{conj:1}.
\begin{defn}
We definite the $p$-energy of $\mu$ to be
\[\cE_p (\mu) := \left(\int \int \frac{1}{|z-w|^p} \, d\mu (z) \, d\mu (w) \right)^{1/p}.\]
Since in the sequel we will only be using the $1$-energy, we will write $\cE$ for $\cE_1.$
\end{defn}
By Fubini's Theorem, when $\mu$ has finite 1-energy, the function
$V_\mu$ defined by
$$V_\mu (z) := \int \frac{1}{z-w} \, d\mu (w)$$
is well defined and in $L^1 (\mu)$.

\begin{thm} \label{th:1}
Suppose $\mu$ has finite 1-energy and that
\begin{equation} \label{eq:nonzero}
\mu \left \{ z : V_\mu (z) = 0 \right \} = 0 \, .
\end{equation}
Then $\zeros (f_n')$ converges in distribution to $\mu$ as
$n \to \infty$.
\end{thm}

A natural set of examples of $\mu$ with finite $1$-energy is provided by the following observation:
\begin{obs}
\label{hausdorff}
Suppose $\Omega \subset \mathbb{C}$ has Hausdorff dimension greater than one, and $\mu$ is in the measure class of the Hausdorff measure on $\Omega.$ Then $\mu$ has finite $1$-energy.
\end{obs}
\noindent{\sc Proof:} 
This is essentially the content of \cite{falcofrac}[Theorem 4.13(b)].
$\Cox$

In particular,  if $\mu$ is uniform in an open subset (with compact 
closure) of $\mathbb{C}$, its $1$-energy is finite.

A natural special case to which Theorem~\ref{th:1} does not apply
is when $\mu$ is uniform on the unit circle; here the 1-energy is
just barely infinite.  
\begin{thm} \label{th:circle}
If $\mu$ is uniform on the unit circle then $\zeros (f_n)$
converges to the unit circle in probability.
\end{thm}

This result is somewhat weak because we do not prove $\zeros (f_n)$
has a limit in distribution, only that all subsequential limits are
supported on the unit circle.
By the Gauss-Lucas Theorem, all roots of $f_n$ have modulus
less than~1, so the convergence to $\mu$ is from the
inside.  Weak convergence to $\mu$ implies that only $o(n)$
points can be at distance $\Theta (1)$ inside the cirle;
the number of such points turns out to be $\Theta (1)$.
Indeed quite a bit can be said about the small outliers.
For $0 < \rho < 1$, define $B_\rho := \{ z : |z| \leq \rho \}$.
The following result, which implies Theorem~\ref{th:circle},
is based on a very pretty result of Peres and 
Virag~\cite[Theorems~1~and~2]{peres-virag} which
we will quote in due course.

\begin{thm} \label{th:pointproc}
For any $\rho \in (0,1)$, as $n \to \infty$, the set
$\zeros (g_n) \cap B_\rho$ of zeros of $g_n$ on $B_\rho$
converges in distribution to a determinantal point process
on $B_\rho$ with the so-called Bergmann kernel 
$\pi^{-1} (1 - z_i \overline{z_j})^2$.
The number $N(\rho)$ of zeros is distributed as the sum of
independent Bernoullis with means $\rho^{2k}$, $1 \leq k < \infty$.
\end{thm}

\subsection{Distance functions on the space of probability measures}
\label{distsec}
If $\mu$ and $\nu$ are probability measures on a separable metric
space $S$, then the \Em{Prohorov}\footnote{Also known as the Prokhorov and the L\'evy-Pro(k)horov distance} distance $|\mu - \nu|_P$ is
defined to be the least $\ee$ such that for every set $A$,
$\mu (A) \leq \nu (A^\ee) + \ee$ and $\nu (A) \leq \mu (A^\ee) + \ee$.
Here, $A^\ee$ is the set of all points within distance $\ee$ of
some point of $A$.  The Prohorov metric metrizes convergence in
distribution.  We view collections of points in $\C$ (e.g., the 
zeros of $f_n$) as probability measures on $\C$, therefore the
Prohorov metric serves to metrize convergence of zero sets.
The space of probability measures on $S$, denoted ${\cal P} (S)$,
is itself a separable metric space, therefore one can define the 
Prohorov metric on ${\cal P} (S)$, and this metrizes convergence
of laws of random zero sets.  

The Ky Fan metric on random variables
on a fixed probability space will be of some use as well.  Defined
by $K(X,Y) = \inf \{ \ee :  \P (d(X,Y) > \ee) < \ee \}$, this
metrizes convergence in probability.  The two metrics are related
(this is Strassen's Theorem):
\begin{equation} \label{eq:prohorov=ky}
|\mu - \nu|_P = \inf \{ K(X,Y) : X \sim \mu , \, Y \sim \nu \} \, .
\end{equation}
A good reference for the facts mentioned above is \cite{hofinger-proh}\footnote{available on line at \url{http://epub.oeaw.ac.at/0xc1aa500d_0x00239061.pdf}}
We will make use of Rouch\'e's Theorem.  There are a number of 
formulations, of which the most elementary is probably the
following statement proved as Theorem~10.10 in~\cite{bak-newman}.
\begin{thm}[Rouch\'e]
If $f$ and $g$ are analytic on a topological disk, $B$, and $|g| < |f|$
on $\partial B$, then $f$ and $f+g$ have the same number of 
zeros on $B$.  
$\Cox$
\end{thm}

\setcounter{equation}{0}
\section{Proof of Theorem~\protect{\ref{th:1}}} \label{sec:results}

We begin by stating some lemmas.  The first is nearly a triviality.
\begin{lem} \label{lem:distance}
Suppose $\mu$ has finite 1-energy.  Then 
\begin{enumerate} \romenumi
\item $$t \cdot \P \left ( |X_0 - X_1| \leq \frac{1}{t} \right ) \to 0 \, .$$
\item for any $C > 0$,
$$\P \left ( \min_{1 \leq j \leq n} |X_j - X_{n+1}| 
   \leq \frac{C}{n} \right ) \to 0 \, ;$$
\end{enumerate}
\end{lem}

\noindent{\sc Proof:} For part~$(i)$ observe that 
$\limsup t \cdot \P ( |X_0 - X_1| \leq 1/t) \leq 
2 \limsup 2^j \cdot \P \left ( |X_0 - X_1| \leq 2^{-j} \right )$
as $t$ goes over reals and $j$ goes over integers.  We then have
\begin{eqnarray*}
\infty & > & \cE (\mu) \\[1ex]
& = & \E 1 / |X_0 - X_1|  \\[1ex]
& \geq & \half \E \sum_{j \in \Z} 2^j \one_{|X_0 - X_1| \leq 2^{-j}} \\[1ex]
& = & \half \sum_j 2^j \P \left ( |X_0 - X_1| \leq 2^{-j} \right )
\end{eqnarray*}
and from the finiteness of the last sum it follows that the summand 
goes to zero.  Part~$(ii)$ follows from part~$(i)$ upon observing,
by symmetry, that 
$$\P \left ( \min_{1 \leq j \leq n} |X_j - X_{n+1}| 
   \leq \frac{C}{n} \right ) 
\leq n \, \P \left ( |X_0 - X_1| \leq \frac{C}{n} \right ) \, .$$
$\Cox$

Define the $n^{th}$ empirical potential function $V_{\mu , n}$ by
$$V_{\mu , n} (z) := \frac{1}{n} \sum_{j=1}^n \frac{1}{z - X_j}$$
which is also the integral in $w$ of $1/(z-w)$ against the measure
$\zeros (f_n)$.  
Our next lemma bounds $V_{\mu , n}' (z)$ on the disk 
$B := B_{C/n} (X_{n+1})$.  
\begin{lem} \label{lem:second}
For all $\ee > 0$, 
$$\P \left ( \sup_{z \in B} |V_{\mu , n}' (z)| 
   \geq \ee n \right ) \to 0$$
as $n \to \infty$.
\end{lem}

\noindent{\sc Proof:} Let $G_n$ denote the event that 
$\min_{1 \leq j \leq n} |X_j - X_{n+1}| > 2 C/n$.  Let
$S_n := \sup_{z \in B} |V_{\mu , n}' (z)|$.  We will show that 
\begin{equation} \label{eq:phi}
\E S_n \one_{G_n} = o(n)
\end{equation}
$n \to \infty$.  By Markov's inequality, this implies that
$\P (S_n \one_{G_n} \geq \ee n) \to 0$ for all $\ee > 0$ as
$n \to \infty$.  By part~$(ii)$ of Lemma~\ref{lem:distance} 
we know that $\P (G_n) \to 1$, which then establishes that
$\P (S_n \geq \ee n) \to 0$, proving the lemma.

In order to show~\eqref{eq:phi} we begin with
$$|V_{\mu , n}' (z)| = \left | \frac{1}{n} \sum_{j=1}^n
   \frac{-1}{(z - X_j)^2} \right |  
   \leq \frac{1}{n} \sum_{j=1}^n \frac{1}{|z - X_j|^2} \, .$$
Therefore, 
$$S_n \one_{G_n} \leq \frac{1}{n} \sum_{j=1}^n 
   \frac{1}{(|X_{n+1} - X_j| - C/n)^2} \leq 
   \frac{1}{n} \sum_{j=1}^n \frac{4}{|X_{n+1} - X_j|^2}$$
and by symmetry,
$$\E S_n \one_{G_n} \leq 
   4 \E |X_0 - X_j|^{-2} \one_{|X_0 - X_j| \geq 2C/n} \leq 
   4 \E \left ( \frac{1}{|X_0 - X_j|} \wedge \frac{2C}{n} \right )^2 
   \, .$$
We may then integrate by parts, obtaining
$$\E S_n \one_{G_n} \leq \int_0^{2 C/n} 8t \, \P \left (
   \frac{1}{|X_0 - X_1|} > t \right ) \, dt \, .$$
The integrand goes to zero as $n \to \infty$ by part~$(i)$ 
of Lemma~\ref{lem:distance}.  It follows that the integral
is $o(n)$, proving the lemma.  $\Cox$

Define the lower modulus of $V$ to distance $C/n$ by
$$\vbar^C_n (z) := \inf_{w : |w-z| \leq C/n} 
   \left | V_{\mu , n} (w) \right | \; .$$
This depends on the argument $\mu$ as well as $C$ and $n$ but we 
omit this from the notation.  

\begin{lem} \label{lem:C}
Assume $\mu$ has finite 1-energy.  Then as $n \to \infty$,
the random variable $\vbar^C_n (X_{n+1})$ converges in probability,
and hence in distribution, to $|V_\mu (X_{n+1})|$.
\end{lem}


In the sequel we will need the Glivenko-Cantelli 
Theorem~\cite[Theorem~1.7.4]{durrett}.  Let $X_1\dotsc, X_n, \dotsc$ 
be independent, identitically distributed random variables in 
$\mathbb{R}$ with common cumulative distribution function $F$.
The \emph{empirical distribution function} $F_n$ for 
$X_1, \dotsc, X_n$ is defined by 
$$F_n(x) = \frac{1}{n} \sum_{i=1}^n I_{(-\infty, x]}(X_i),$$
where $I_C$ is the indicator function of the set $C.$ For every fixed $x$, $F_n(x)$ is a sequence of random variables, which converges to $F(x)$ almost surely by the strong law of large numbers. Glivenko-Cantelli Theorem strengthen this by proving \emph{uniform convergence} of $F_n$ to $F.$
\begin{thm}[Glivenko-Cantelli]
\label{glican}
\[
\|F_n -F\|_\infty = \sup_{x\in \mathbb{R}} \left|F_n(x) - F(x)\right| \longrightarrow 0\quad \mbox{almost surely.}
\]
\end{thm}
The following Corollary is immediate:
\begin{cor}
\label{glicancor}
Let $f$ be a bounded continuous function on $\mathbb{R}.$ Then
\[
\lim_{n\rightarrow \infty} \int_{\mathbb{R}} f dF_n = \int_{\mathbb{R}} f d F, \quad \mbox{almost surely}.
\]
\end{cor}
Another immediate Corollary is:
\begin{cor}
\label{glicanpro}
With notation as in the statement of Theorem \ref{glican}, the Prohorov distance between $F_n$ and $F$ converges to zero almost surely.
\end{cor}
\noindent{\sc Proof of Lemma \ref{lem:C}:} It is equivalent to show that 
$\vbar^C_n - | V_\mu (X_{n+1}) | \to 0$ in probability, 
for which it sufficient to show
\begin{equation} \label{eq:to show 1}
\sup_{u \in B} \left | V_{\mu,n} (u) - V_\mu (X_{n+1}) \right | \to 0
\end{equation}  
in probability.  This will be shown by proving the following
two statements:
\begin{eqnarray}
\sup_{u \in B} \left | V_{\mu,n} (u) - V_{\mu,n} (X_{n+1}) 
   \right | & \to & 0 
   \mbox{ in probability} \; ; \label{eq:to show 2} \\[2ex]
\left | V_{\mu,n} (X_{n+1}) - V_\mu (X_{n+1}) \right | & \to & 0
   \mbox{ in probability} \label{eq:to show 3} \, .
\end{eqnarray}
The left-hand side of~\eqref{eq:to show 2} is bounded above by
$(C / n) \sup_{u \in B} |V_{\mu,n}' (u)|$.  By
Lemma~\ref{lem:second}, for any $\ee > 0$, the probability of this
exceeding $C \ee$ goes to zero as $n \to \infty$.  This 
establishes~\eqref{eq:to show 2}.

For~\eqref{eq:to show 3} we observe, using Dominated Convergence, that
under the finite 1-energy condition, 
$$\cE^K (\mu) := \int \int \frac{1}{|z-w|} \one_{|z-w|^{-1} \geq K}
    \, d\mu (z) \, d\mu (w) \to 0$$
as $K \to \infty$.  Define $\phi^{K,z}$ by
$$\phi^{K,z} (w) = \frac{1}{z-w} \frac{|z-w|}{\max \{ |z-w| , 1/K \}}$$
in other words, it agrees with $1/(z-w)$ except that we 
multiply by a nonegative real so as to truncate the magnitude
at $K$.  We observe for later use that 
$$\left | \phi^{K,z} (w) - \frac{1}{|z-w|}  \right | \leq 
   \frac{1}{|z-w|} \one_{|z-w|^{-1} \geq K}$$
so that 
\begin{equation} \label{eq:cEK}
\int \int \left | \phi^{K,z} (w) - \frac{1}{|z-w|} \right |
   \, d \mu (z) \, d\mu (w) \leq \cE^K (\mu)  \to 0 \, .
\end{equation}

We now introduce the truncated potential and truncated empirical
potential with respect to $\phi^{K,z}$:
\begin{eqnarray*}
V_\mu^K (z) & := & \int \phi^{K,z} (w) \, d\mu (w) \\[1ex]
V_{\mu,n}^K (z) & := & \int \phi^{K,z} (w) \, d\zeros (f_n) (w) \, .
\end{eqnarray*}
We claim that 
\begin{equation} \label{eq:E}
\E \left | V^K_\mu (X_{n+1}) - V_\mu (X_{n+1}) \right | \leq \cE^K (\mu) \, .
\end{equation}
Indeed, 
$$V_\mu (X_{n+1}) - V^K_\mu (X_{n+1}) = 
   \int \left ( \frac{1}{z-X_{n+1}} - \phi^{K,z} (X_{n+1}) \right ) 
   \, d\mu (z) \, $$
so taking an absolute value inside the integral, then integrating
against the law of $X_{n+1}$ and using~\eqref{eq:cEK} 
proves~\eqref{eq:E}.  The empirical distribution $V_{\mu , n}$
has mean $\mu$ and is independent of $X_{n+1}$, therefore
the same argument proves
\begin{equation} \label{eq:En}
\E \left | V^K_{\mu,n} (X_{n+1}) - V_{\mu,n} (X_{n+1}) \right | 
   \leq \cE^K (\mu) 
\end{equation}
independent of the value of $n$.  

We now have two thirds of what we need for the triangle inequality.  
That is, to show~\eqref{eq:to show 3} we will show that the following
three expressions may all be made smaller than $\ee$ with probability
$1 - \ee$.  
\begin{eqnarray*}
&& V_{\mu,n} (X_{n+1}) - V^K_{\mu,n} (X_{n+1}) \\[2ex]
&& V^K_{\mu,n} (X_{n+1}) - V^K_\mu (X_{n+1}) \\[2ex]
&& V^K_{\mu} (X_{n+1}) - V_\mu (X_{n+1}) 
\end{eqnarray*}
Choosing $K$ large enough so that $\cE^K (\mu) < \ee^2$, this
follws for the third of these follows by~\eqref{eq:E} and for
the first of these by~\eqref{eq:En}.  Fixing this value of
$K$, we turn to the middle expression.  The function $\phi^{K,z}$
is bounded and continuous.  By the Corollary \ref{glicancor} to the Glivenko-Cantelli Theorem \ref{glican}, 
the empirical law $\zeros (f_n)$ converges weakly to $\mu$, meaning
that the integral of the any bounded continuous function $\phi$
against $\zeros (f_n)$ converges in probability to the integral
of $\phi$ against $\mu$.  Setting $\phi := \phi^{K,z}$ and $z :=
X_{n+1}$ proves that $V^K_{\mu,n} (X_{n+1}) - V^K_\mu (X_{n+1})$ 
goes to zero in probability, establishing the middle statement
(it is in fact true conditionally on $X_{n+1}$) and concluding 
the proof.
$\Cox$

\noindent{\sc Proof of Theorem}~\ref{th:1}:
Suppose that $\vbar_n^C (X_{n+1}) > 1/C$.  Then for all $w$ with 
$|w - X_{n+1}| \leq C/n$, we have
$$f_n' (w) = \sum_{j=1}^n \frac{1}{w - X_j} = n V_{\mu , n} (w)
   \geq \frac{n}{C}$$
and hence
$$\left | f_n' (w) \right | = n \left | V_{\mu,n} (w) \right |
   \geq n \vbar_n^C (X_{n+1}) \geq \frac{n}{C} \, .$$
To apply Rouch\'e's Theorem to the functions $1/f_n'$ and 
$z - X_{n+1}$ on the disk $B := B_{C/n} (X_{n+1})$ we note
that $|1/f_n'| < C / n = |z-X_{n+1}|$ on $\partial B$
and hence that the sum has precisely one zero in $B$,
call it $a_{n+1}$.  Taking reciprocals we see that $a_{n+1}$ 
is also the unique value in $z \in B$ for which $f_n' (z) = 
-1/(z - X_{n+1})$.  But $f_n' (z) + 1/(z - X_{n+1}) = f_{n+1}' (z)$,
whence $f_{n+1}'$ has the unique zero $a_{n+1}$ on $B$.

Now fix any $\delta > 0$.  Using the hypothesis that 
$\mu \{ z : V_\mu (z) = 0 \} = 0$, we pick a $C > 0$
such that $\P (|V_\mu (X_{n+1})| \leq 2/C) \leq \delta/2$.
By Lemma~\ref{lem:C}, there is an $n_0$ such that for all
$n \geq n_0$,
$$\P \left ( \vbar^C (X_{n+1}) \leq \frac{1}{C} \right ) 
   \leq \delta \, .$$
It follows that the probability
that $f_{n+1}'$ has a unique zero $a_{n+1}$ in $B$ is at
least $1 - \delta$ for $n \geq n_0$.  By symmetry, we see 
that for each $j$, the probability is also at least $1-\delta$
that $f_{n+1}'$ has a unique zero, call it $a_j$, in the ball 
of radius $C/n$ centered at $X_j$; equivalently, the expected 
number of $j \leq n+1$ for which there is not a unique zero of 
$f_{n+1}'$ in $B_{C/n} (X_j)$ is at most $\delta n$ for $n \geq n_0$.

Define $x_j$ to equal $a_j$ if $f_{n+1}'$ has a unique root
in $B_{C/n} (X_j)$ and the minimum distance from $X_j$ to
any $X_i$ with $i \leq n+1$ and $i \neq j$ is at least $2C/n$.
By convention, we define $x_j$ to be the symbol $\Delta$ if
either of these conditions fails.  The values $x_j$ other 
than $\Delta$ are distinct roots of $f_{n+1}'$ and each such
value is within distance $C/n$ of a different root of $f_{n+1}$.  
Using part~$(ii)$ of Lemma~\ref{lem:distance} we see that the 
expected number of $j$ for which $x_j = \Delta$ is $o(n)$.  It
follows that $\P (|\zeros (f_{n+1}) - \zeros (f_{n+1}')|_P 
\geq 2 \delta) \to 0$ as $n \to \infty$.
But also the Prohorov distance between $\zeros (f_{n+1})$ 
and $\mu$ converges to zero by Corollary \ref{glicanpro}.
The Prohorov distance metrizes convergence in distribution
and $\delta > 0$ was arbitrary, so the theorem is proved.
$\Cox$

\setcounter{equation}{0}
\section{Proof of remaining theorems}

\noindent{\sc Proof of Theorem}~\ref{th:pointproc}:
Let $\G := \sum_{j=0}^\infty Y_j z^j$ denote the standard 
complex Gaussian power series where $\{ Y_j (\omega) \}$ are 
IID standard complex normals.  The results we require 
from~\cite{peres-virag} are as follows.
\begin{pr}[\protect{\cite{peres-virag}}] \label{pr:PV}
The set of zeros of $\G$ in the unit disk is a determinantal 
point process with joint intensities 
$$p(z_1 , \ldots , z_n) = \pi^{-n} \det \left [ 
   \frac{1}{(1 - z_i \overline{z}_j)^2} \right ] \; .$$
The number $N(\rho)$ of zeros of $\G$ on $B_\rho$ is distributed 
as the sum of independent Bernoullis with means $\rho^{2k}$, 
$1 \leq k < \infty$. 
$\Cox$
\end{pr}

To use these results we broaden them to random series whose
coefficients are nearly IID Gaussian.
\begin{lem} \label{lem:series}
Let $\{ g_n := \sum_{r=0}^\infty a_{nr} z^r \}$ be a sequence
of power series.  Suppose
\begin{enumerate} \romenumi
\item for each $k$, the $k$-tuple $(a_{n,1} , \ldots , 
a_{n,k})$ converges weakly as $n \to \infty$ to a $k$-tuple of IID 
standard complex normals;  
\item $\E |a_{nr}| \leq 1$ for all $n$ and $r$. 
\end{enumerate}
Then on each disk $B_\rho$, the set $\zeros (g_i) \cap B_\rho$ 
converges weakly to $\zeros (\G) \cap \rho$.  
\end{lem}

\noindent{\sc Proof:} Throughout the proof we fix $\rho \in (0,1)$
and denote $B := B_\rho$.  Suppose an analytic function $h$ has no 
zeros on $\partial B$.  Denote by $||g-h||_B$ the sup norm on
functions restricted to $B$.  Note that if $h_n \to h$ uniformly on 
$B$ then $\zeros (h_n) \cap B \to \zeros (h) \cap B$ in the 
weak topology on probability measures on $B$, provided that $h$
has no zero on $\partial B$.  We apply this with
$h = \G := \sum_{j=0}^\infty Y_j z^j$ where $\{ Y_j (\omega) \}$ are IID
standard complex normals.  For almost every $\omega$, $h(\omega)$
has no zeros on $\partial B$.  Hence given $\ee > 0$ there is almost
surely a $\delta (\omega) > 0$ such that $||g - \G||_B < \delta$ implies
$|\zeros (g) - \zeros (\G)|_P < \ee$.  Pick $\delta_0 (\ee)$
small enough so that $\P (\delta(\omega) \leq \delta_0) < \ee/3$;
thus $||g-\G||_B < \delta_0$ implies $|\zeros (g) - \zeros (\G)| < \ee$
for all $\G$ outside a set of measure at most $\ee / 3$.  

By hypothesis~$(ii)$, 
$$\E \left | \sum_{r=k+1}^\infty a_{nr} z^r \right | 
   \leq \frac{\rho^{k+1}}{1 - \rho} \, .$$
Thus, given $\ee > 0$, once $k$ is large enough so that 
$\rho^{k+1}/(1-\rho) < \ee \delta_0 (\ee) / 6$, we see that 
$$\P \left ( \left | \sum_{r=k+1}^\infty a_{nr} z^r \right | 
   \geq \frac{\delta_0 (\ee)}{2} \right ) \leq \frac{\ee}{3} \, .$$
For such a $k (\ee)$ also $|\sum_{r=k+1}^\infty Y_r z^r| \leq \ee / 3$.
By hypothesis~$(i)$, given $\ee > 0$ and the corresponding 
$\delta (\ee)$ and $k(\ee)$, we may choose $n_0$ such that $n \geq n_0$
implies that the law of $(a_{n1} , \ldots , a_{nk})$ is within
$\min \{ \ee / 3 , \delta_0 (\ee) / (2k) \}$ of the product
of $k$ IID standard complex normals in the Prohorov metric.
By the equivalence of the Prohorov metric to the minimal
Ky Fan metric, there is a pair of random variables $\gt$ and $\htt$
such that $\gt \sim g_n$ and $\htt \sim \G$ and, except on a set of
of measure $\ee/3$, each of the first $k$
coefficients of $\gt$ is within $\delta_0/(2k)$ of the corresponding 
coefficient of $\G$.  By the choice of $k (\ee)$, we then have
$$\P (||\gt - \htt||_B \geq \delta_0) \leq \frac{2 \ee}{3} \, .$$ 
By the choice of $\delta_0$, this implies that
$$\P (|\zeros (\gt) - \zeros (\htt)|_P \geq \ee) < \ee \, .$$
Because $\gt \sim g_n$ and $\htt \sim \G$, we see that
the law of $\zeros (g_n) \cap B$ and the law of $\zeros (\G) \cap B$ 
are within $\ee$ in the Prohorov metric on laws on measures.  Because
$\ee > 0$ was arbitrary, we see that the law of $\zeros (g_n) \cap B$ 
converges to the law of $\zeros (\G) \cap B$.
$\Cox$

\noindent{\sc Proof of Theorem}~\ref{th:pointproc}:
Let $\rho < 1$ be fixed for the duration of this argument and denote
$B := B_\rho$.  Let
$$g_n (z) := \frac{f_n' (z)}{f(z)} = \sum_{j=1}^n \frac{1}{z - X_j} 
   \, .$$ 
Because $|X_j| = 1$, the rational function $1/(z - X_j) = 
-X_j^{-1} / (1 - X_j^{-1} z)$ is analytic on the open unit disk 
and represented there by the power series $- \sum_{r=0}^\infty 
X_j^{-r-1} z^r$.  It follows that $- g_n / \sqrt{n}$ is analytic 
on the open unit disk and represented there by the power series 
$- g_n (z) / \sqrt{n} = \sum_{r=0}^\infty a_{nr} z^r$ where
$$a_{nr} = n^{-1/2} \sum_{j=1}^n X_j^{-r-1} \, .$$
The function $- g_n / \sqrt{n}$ has the same zeros on $B$ as does 
$f_n'$, the normalization by $- 1 / \sqrt{n}$ being inserted as
a convenience for what is about to come.

We will apply Lemma~\ref{lem:series} to the sequence $\{ g_n \}$.
The coefficients $a_{nj}$ are normalized power sums of the
variables $\{ X_j \}$.  For each $r \geq 0$ and each $j$, the 
variable $X_j^{-r-1}$ is uniformly distributed on the unit circle.
It follows that $\E a_{nr} = 0$ and that $\E a_{nr} \overline{a_{nr}} 
= n^{-1} \sum_{ij} X_i^{-r-1} \overline{X_j}^{-r-1} 
= n^{-1} \sum_{ij} \delta_{ij} = 1$.  In particular, 
$\E |a_{nr}| \leq (\E |a_{nr}|^2)^{1/2} = 1$, satisfying
the second hypothesis of Lemma~\ref{lem:series}.  For
the first hypothesis, fix $k$, let $\theta_j = \Arg (X_j)$, and
let $\vv^{(j)}$ denote the $(2k)$-vector $(\cos (\theta_j), - \sin (\theta_j),
\cos (2 \theta_j) , - \sin (2 \theta_j) , \ldots , \cos (k \theta_j) , 
- \sin (k \theta_j))$; in other words, $\vv^{(j)}$ is the complex
$k$-vector $(X_j^{-1} , X_j^{-2} , \ldots , X_j^{-k})$ viewed
as a real $(2k)$-vector.  For each $1 \leq s, t \leq 2k$ we have
$\E \vv^{(j)}_s \vv^{(j)}_t = (1/2) \delta_{ij}$.  Also the
vectors $\{ \vv^{(j)} \}$ are independent as $j$ varies.  It
follows from the multivariate central limit theorem
(see, e.g.,~\cite[Theorem~2.9.6]{durrett}) that $\uu^{(n)} 
:= n^{-1/2} \sum_{j=1}^n \vv^{(j)}$ converges to $1/\sqrt{2}$ 
times a standard $(2k)$-variate normal.  
For $1 \leq r \leq k$, the coefficient $a_{nr}$ is equal to
$\uu^{(n)}_{2r-1} + i \uu^{(n)}_{2r}$.  Thus $\{ a_{nr} : 
1 \leq r \leq k \}$ converges in distribution as $n \to \infty$
to a $k$-tuple of IID standard complex normals.  The hypotheses of
Lemma~\ref{lem:series} being verified, the theorem now follows
from Proposition~\ref{pr:PV}.
$\Cox$

\bibliographystyle{alpha}
\bibliography{RP}

\def\cdprime{$''$}
\begin{thebibliography}{GRR69}

\bibitem[Azi85]{MR788579}
Abdul Aziz.
\newblock On the zeros of a polynomial and its derivative.
\newblock {\em Bull. Austral. Math. Soc.}, 31(2):245--255, 1985.

\bibitem[BN82]{bak-newman}
J.~Bak and D.~Newman.
\newblock {\em Complex Analysis}.
\newblock undergraduate Texts in Mathematics. Springer-Verlag, Berlin, 1982.

\bibitem[Bra31]{MR1506859}
Hubert~E. Bray.
\newblock On the {Z}eros of a {P}olynomial and of {I}ts {D}erivative.
\newblock {\em Amer. J. Math.}, 53(4):864--872, 1931.

\bibitem[{\'C}M04]{MR2063118}
Branko {\'C}urgus and Vania Mascioni.
\newblock A contraction of the {L}ucas polygon.
\newblock {\em Proc. Amer. Math. Soc.}, 132(10):2973--2981 (electronic), 2004.

\bibitem[dB46]{MR0019157}
N.~G. de~Bruijn.
\newblock On the zeros of a polynomial and of its derivative.
\newblock {\em Nederl. Akad. Wetensch., Proc.}, 49:1037--1044 = Indagationes
  Math. 8, 635--642 (1946), 1946.

\bibitem[dBS47]{MR0021148}
N.~G. de~Bruijn and T.~A. Springer.
\newblock On the zeros of a polynomial and of its derivative. {II}.
\newblock {\em Nederl. Akad. Wetensch., Proc.}, 50:264--270=Indagationes Math.
  9, 458--464 (1947), 1947.

\bibitem[Dim98]{MR1452801}
Dimitar~K. Dimitrov.
\newblock A refinement of the {G}auss-{L}ucas theorem.
\newblock {\em Proc. Amer. Math. Soc.}, 126(7):2065--2070, 1998.

\bibitem[Dro89]{MR1063061}
Janusz Dronka.
\newblock On the zeros of a polynomial and its derivative.
\newblock {\em Zeszyty Nauk. Politech. Rzeszowskiej Mat. Fiz.}, (9):33--36,
  1989.

\bibitem[Dur04]{durrett}
R.~Durrett.
\newblock {\em Probability: Theory and Examples}.
\newblock Duxbury Press, Belmont, CA, third edition, 2004.

\bibitem[Fal03]{falcofrac}
Kenneth Falconer.
\newblock {\em Fractal geometry}.
\newblock John Wiley \& Sons Inc., Hoboken, NJ, second edition, 2003.
\newblock Mathematical foundations and applications.

\bibitem[GRR69]{MR0239062}
A.~W. Goodman, Q.~I. Rahman, and J.~S. Ratti.
\newblock On the zeros of a polynomial and its derivative.
\newblock {\em Proc. Amer. Math. Soc.}, 21:273--274, 1969.

\bibitem[Hof06]{hofinger-proh}
Andreas Hofinger.
\newblock The metrics of {P}rokhorov and {K}y {F}an for assessing uncertainty
  in inverse problems.
\newblock {\em \"Osterreich. Akad. Wiss. Math.-Natur. Kl. Sitzungsber. II},
  215:107--125 (2007), 2006.

\bibitem[Joy69]{MR0239061}
Andr{\'e} Joyal.
\newblock On the zeros of a polynomial and its derivative.
\newblock {\em J. Math. Anal. Appl.}, 26:315--317, 1969.

\bibitem[Mah61]{MR0133437}
K.~Mahler.
\newblock On the zeros of the derivative of a polynomial.
\newblock {\em Proc. Roy. Soc. Ser. A}, 264:145--154, 1961.

\bibitem[Mal05]{MR2159699}
S.~M. Malamud.
\newblock Inverse spectral problem for normal matrices and the {G}auss-{L}ucas
  theorem.
\newblock {\em Trans. Amer. Math. Soc.}, 357(10):4043--4064 (electronic), 2005.

\bibitem[Mar49]{marden1949}
M.~Marden.
\newblock {\em Geometry of Polynomials}, volume~3 of {\em Mathematical Surveys
  and Monographs}.
\newblock AMS, 1949.

\bibitem[Mar83]{MR700266}
Morris Marden.
\newblock Conjectures on the critical points of a polynomial.
\newblock {\em Amer. Math. Monthly}, 90(4):267--276, 1983.

\bibitem[Paw98]{MR1473453}
Piotr Pawlowski.
\newblock On the zeros of a polynomial and its derivatives.
\newblock {\em Trans. Amer. Math. Soc.}, 350(11):4461--4472, 1998.

\bibitem[PV05]{peres-virag}
Y.~Peres and B.~Virag.
\newblock Zeros of the i.i.d. {G}aussian power series: a conformally invariant
  determinantal process.
\newblock {\em Acta Math.}, 194:1--35, 2005.

\bibitem[Rah72]{MR0308374}
Q.~I. Rahman.
\newblock On the zeros of a polynomial and its derivative.
\newblock {\em Pacific J. Math.}, 41:525--528, 1972.

\bibitem[Sen01]{MR1856819}
Bl. Sendov.
\newblock Hausdorff geometry of polynomials.
\newblock {\em East J. Approx.}, 7(2):123--178, 2001.

\bibitem[Sen10]{MR2722582}
Blagovest Sendov.
\newblock New conjectures in the {H}ausdorff geometry of polynomials.
\newblock {\em East J. Approx.}, 16(2):179--192, 2010.

\bibitem[Sto96]{MR1400355}
{\`E}.~A. Storozhenko.
\newblock On a problem of {M}ahler on the zeros of a polynomial and its
  derivative.
\newblock {\em Mat. Sb.}, 187(5):111--120, 1996.

\bibitem[Tar86]{MR874708}
Q.~M. Tariq.
\newblock On the zeros of a polynomial and its derivative. {II}.
\newblock {\em J. Univ. Kuwait Sci.}, 13(2):151--156, 1986.

\end{thebibliography}

\end{document}